\documentstyle[11pt]{article}
\begin{document}
\def\baselinestretch{1.2}
\newcommand{\cA}{{\cal A}}
\newcommand{\cB}{{\cal B}}
\newcommand{\cC}{{\cal C}}
\newcommand{\cD}{{\cal D}}
\newcommand{\cE}{{\cal E}}
\newcommand{\cF}{{\cal F}}
\newcommand{\cG}{{\cal G}}
\newcommand{\cH}{{\cal H}}
\newcommand{\cI}{{\cal I}}
\newcommand{\cJ}{{\cal J}}
\newcommand{\cK}{{\cal K}}
\newcommand{\cL}{{\cal L}}
\newcommand{\cM}{{\cal M}}
\newcommand{\cN}{{\cal N}}
\newcommand{\cO}{{\cal O}}
\newcommand{\cP}{{\cal P}}
\newcommand{\cQ}{{\cal Q}}
\newcommand{\cR}{{\cal R}}
\newcommand{\cS}{{\cal S}}
\newcommand{\cT}{{\cal T}}
\newcommand{\cU}{{\cal U}}
\newcommand{\cV}{{\cal V}}
\newcommand{\cW}{{\cal W}}
\newcommand{\cX}{{\cal X}}
\newcommand{\cY}{{\cal Y}}
\newcommand{\cZ}{{\cal Z}}

\newcommand{\ST}{{-\!\!\!-\!\!\!-\!\!\!-\!\!\!-\!\!\!
-\!\!\!\!\!\!-\!\!\!\!\!\!-\!\!\!\!\!\!-\!\!\!\!\!\!-\!\!\!\!\!\!\longrightarrow}}
\newcommand{\DT}{{\begin{array}{ll}
\!\!\!\!\mid & {}  \\[-2mm]
\!\!\!\!\mid & {} \\ [-2mm]
\!\!\!\!\mid & {}  \\[-2mm]
\!\!\!\!\mid & {}   \\[-2mm]
\!\!\!\!\!\vee & {}
\end{array}
}}
\begin{center}
{\Large{\bf Compactified Jacobian of some Unibranch curves}}
Jyotsna Gokhale\\
Mehta Research Institute\\
10 Kasturba Gandhi Marg\\
Allahabad 211002\\
India.\\
E-mail : jgokhale@mri.ernet.in
\end{center}
\noindent{\underline{\bf Introduction}}

Let $C$ be an integral curve with one unibranch singularity of
embedding dimension greater than two.  It is known that
$ \ (Pic^d C)^{=} \ $,
the moduli of torsion free sheaves of rank one and degree $d$,
is irreducible if and only if $C$  is locally planar ([D], [AIK], [R], [K]).
It is enough to consider a curve with one singularity [SK].

In this work
\footnote{
This work was done during '84 - '86  at Brandeis University as part of my
doctoral dissertation,and for various reasons remained unpublished so far
along with rest of my dissertation. It was brought to my notice in '92,
when this paper was sent to Mathematische Zeitschrift for publication,
that G.Pfister and J.H.M.Steenbrink have recently ('91) published a paper
similar in content.  We have consequently shortened the paper by giving
the reference of Phister and Steenbrink whenever possible, omitting
repetition of already published proofs.

It is a duty and a pleasure to thank my advisor David Eisenbud for all the
insight and the constructive comments he provided me with, and for the
atmosphere he generates around him.}
we give a basic decomposition of
$ \ (Pic^0 C)^{=} \ $
towards irreducible components for an integral curve with one unibranch
singularity, and give the decomposition into irreducible components
for some types of curves.
$\hfill{\Diamond}$
\newpage
\underline{\bf Notation :}

\begin{tabular}{lll}
$k$                     & {} & an algebraically closed field
                                      of characteristic zero \\
$C$                     & {} & an integral curve over $k$ with
                                 one unibranch singularity $Q$ \\
$C^{\sim}$              & {} & normalisation of $C$ \\
$\pi$                   & {} & the normalisation map \\
$\cO$                   & {} & completion of local ring of $C$ at $Q$ \\
$\cO^{\sim}$            & {} & normalisation of $\cO$ \\
$\cM$                   & {} & maximal ideal of $\cO$ \\
$\cC$                   & {} & the conductor $Hom_{\cO}(\cO^{\sim},\cO)$ of $C$
\\
$\cM_Q$                 & {} & maximal ideal sheaf of $Q$ \\
$\cC_Q$                 & {} & conductor sheaf \\
$\delta$                & {} & $ rk_k (\cO^{\sim} / \cO);$
                                we omit the subscript k normally \\
$\cF$                   & {} & torsion free sheaf of rank one \\
$\cL $                  & {} & a locally free sheaf of rank $1$ \\
$\cW^0_C$               & {} & the canonical dualising sheaf of $C$ \\
$C'$                    & {} & first partial normalization of $C$ such that \\
    {}                  & {} & $ \ rk_k(\cO'/ \cO) \ = \ 1 \ $
                                 and $ \cC' \neq \cC \ $;
                          i.e.,  $ \ \cO' \ = \ \cO + t^{-1}\cC \ $.
\\
$C^i$                   & {} & successive partial normalizations
                                     defined along the same line \\
$\cM^i$                 & {} & maximal ideal of $\cO^i$ \\
$\cC^i$                 & {} & the conductor of $C^i$ \\
$Pic^0 C$               & {} & Picard group of degree $0$ of $C$ \\
$\pi^*$                 & {} & the map $Pic^0 C \longrightarrow Pic^0 C^{\sim}$
                                    of Picard groups corresponding to $\pi$ \\
$K$                     & {} & $ \ ker (\pi^*)$  \\
$(Pic^0 C)^{=}$         & {} & moduli of torsion free sheaves
                                             of degree $0$ and rank one. \\
$\overline{S}$          & {} & Zariski closure of a variety or scheme $S$\\
$(\overline{Pic^0 C})$  & {} & Zariski closure of $Pic^0C$ in $(Pic^0 C)^{=}$\\
$B(\overline{Pic^0 C})$ & {} & boundary of $\overline{Pic^0 C}$, def(0.1) \\
$\chi(\cF)$             & {} &  Euler characteristic of $\cF$ \\
$deg(\cF)$              & {} & $\chi(\cF)-\chi(\cO_C)$ \\
$\Gamma$                & {} & $\Gamma (C) = (k_1, \cdots, k_r)\ {\bf N}$,
                                     the semigroup of orders of generators \\
      {}                & {} & of $\cO$ over $k$, assuming
                                $ 2 < k_1 < k_2 < \cdots < k_r $ \\
$v_0$                   & {} & order of the conductor, i.e.,
                                   $ \cC = t^{v_0}\cO^{\sim} $ \\
$t$                     & {} & any local parameter of
                                 $ \cO^{\sim}. \hfill{\Diamond} $ \\
\end{tabular}
\footnote{$\Gamma^*$ of [PS] corresponds to our $C'.$}
\newpage
\noindent{\underline{\bf  0: Preliminaries}}

We define
$ \ (Pic_{C}^{d})^=(S) \ = \ \{ $
isomorphism classes of torsion free
$ \ \cO_C(S) $
-modules of degree $d$ and rank 1 \}.
$ \ (Pic^{d}_{C})^{=} \ $ is representable and is
represented by a projective scheme
$ \ (Pic^d C)^{=} \ $ \
[D], [AK], which is irreducible [R], [K], iff
$ \ rk(\cM / \cM^2) \leq 2. \ $
We define the boundary of
$ \ \overline{Pic \ C}, \ $ a la [R]:
\medskip\\
{\bf Definition 0.1:}  $ \ \cF \in B(\overline{Pic^0 \ C}) \ $  if
$ \ \cF \not \in Pic^0 C, \ $
there exists an
$ \cO_{C \chi \ Spec \ k[[t]]} $
-module
$ \cG \ $
flat over
$ k[[t]] \ $
such that
$ \ (\cG/t\cG) \approx \cF \ $
and the stalk of $\cG$ at the generic point of
$ \ Spec \ k[[t]] \ $
is a locally free
$\cO_C {\rm -module}. \hfill{\parallel}$
\medskip\\
Let
$ K = \ ker (\pi^*)$,
where
$ \ \pi^* \ :  \ Pic^0 C \longrightarrow Pic^0 C^{\sim} \ $
is the map of respective Picard groups corresponding to the normalization map
$ \ \pi : C^{\sim} \longrightarrow C. $

\noindent{\bf [R,1.2]:}  $ \ \cF \ \in \ \overline{Pic^0 C} \
\Longleftrightarrow \ \exists \ \cL \ \in \ Pic^0 C \ $ such that
$ \ \cF \otimes \cL \ \in \ \overline{K}. \hfill{\parallel}$

For
$ d \leq rk(\cO^{\sim} / \cC), $
the functor
$$
\begin{array}{lll}
E(\cC, d)(S) & = & \{ {\rm isomorphism \ classes \ of}
\ \cO_S {\rm -modules} \ F_S \mid  \\ {} & {} &
\cC \otimes_k \cO_S \subseteq F_S \subseteq \cO^{\sim} \otimes \cO_S,
{\rm and} \ rk(\cO^{\sim} \otimes \cO_S / F_S)  =  d \}
\end{array}
$$
is representable by a projective scheme [R],
also called
$ E(\cC, d); $
we identify
$ E(\cC, d) $
with its image in
$ (Pic^0 C)^{=} $.
In particular, since
$ K \subseteq E(\cC, \delta), \ $
every boundary point ``defines an element of
$ \ E(\cC, \delta)''\ $ [R, Th.2.3(b)].
So it is enough to consider limits of one - parameter families of free
$\cO$-modules
$ \ \{ F_{\beta} \} \ \subseteq \ E(\cC, \delta) \ $
parameterizsed by
$ \ k[[\beta]]. \ $
A free module
$ \ \cG \ \in \ E(\cC, \delta)(k(\beta)) \ $
is necessarily of the kind
$ \ \cG \ = \ \partial_{\beta} \cO, \
\partial_{\beta} \ \in \ [\cO \otimes k(\beta)]^* \ $.
We shall write
$ \ \partial_{\beta}\cO \ \longrightarrow \ F \ $
if $F$ is in the closure of such a family.

A key concept of help in determining
$ B(\overline{Pic^0 C}) $
from [J] is that of the unique filtration of $\cO$
given by
$ \ \cO \supseteq \cM \supseteq \cdots \supseteq I_j \cdots \supseteq \cC \ $
such that
$ \ rk_k(I_j/I_{j+1}) = 1, \ {\rm and} \
ord_t(I_j) < ord_t(I_{j+1}) \ \ \forall \ j \geq 0. $
\footnote{The definition of
$ I_n $
in [PS] is parallel to this.}
Every
$ F \in E(\cC, \delta) \ $
has such a unique filtration.
We define [J'] a functor
$ Filt_{\cC}^{\delta} $
using this concept of filtrations of
$ \cF $
and proved that it is representable by a projective scheme
$ \ Filt(\cC, \delta) \ = \  \{ F \in E(\cC, \delta) \mid
F_j \ \subseteq \ I_j \cO^{\sim} \ \ \forall \ j \geq 0 \}. \ $
Since
$ \overline{K} \subseteq Filt(\cC, \delta) \subseteq E(\cC,\delta) $
this is a major tool in determining
$ B(\overline{Pic^0 C}). $
For brevity we shall use the notation
$ Filt(\cC, \delta) $
in this work.
\smallskip\\
We will close the section with some remarks implicitly assumed and
used in further sections. Since
$ v_0 - \delta $
is used frequently, we will use notation
$ v_0 - \delta = \gamma \ $
for ease of reading.
\medskip\\
\underline{\bf Remarks:}
\begin{enumerate}
\item[(0.2)] Since
$ \# \{ j \in {\bf N} \mid j \not \in \Gamma \} = \delta, $
so
$ \# \{ j \in \Gamma \mid j < v_0 \} = \gamma. $
\item[(0.3)]
$ j \in \Gamma
\Longrightarrow v_0 -1-j  \not \in \Gamma \ \forall j \in {\bf Z} $.
The converse is true iff $C$ is Gorenstein [HK].
\item[(0.4)] (i)
$ \delta - 1 \leq v_0 \leq 2 \delta $. \\
(ii)
$ \# \{ j \in {\bf N} \mid j, v_0 -1-j \not \in \Gamma \} = 2\delta -v_0$.
\item[(0.5)] If $F$ is an $\cO$-module such that
$ \cC \subseteq F \subseteq \cO^{\sim} \ $,
and
$ End(F) $ does not contain
$ t^{-1} \cO \ $,
then
$ \cO \subset F \subseteq \cO^{\sim}, \ $
and
$ v_0 - \delta \leq rk(F / \cC) \leq \delta. $
Further, if
$ rk(F / \cC) = \delta $
then
$ F \approx \cW^{0}_{C,Q} \ $ necessarily. Conversely, for
$ j : v_0 - \delta \leq \ j \leq \delta, $
there exists an $F$ such that
$ \cC \subseteq F \subseteq \cO^{\sim}, \ rk(F / \cC) = j, \ $
and
$ End(F) $
does not contain
$ t^{-1} \cC. $
\item[(0.6)] The following are equivalent: \\
(i) $ v_0 = k_{\gamma-1} + k_1. $ \\
(ii) $ k_i = ik_1 \ \forall \ i \leq \gamma, \ k_i \in \Gamma. $ \\
(iii) $ \ rk(\cM / \cM^2 + t\cC) = 1.$
\item[(0.7)]
$ E(\cC, \delta) \ {\rm is \ irreducible}
\ \Longleftrightarrow \ \overline{K} \ = \ E(\cC, \delta). \ $
\item[(0.8)] If
$ Hom_{\cO_c} (\cF, \cO_c) \ = \ \cF^v \in B(\overline{Pic^0 C}) \ $
then
$ \ F \in \ B(\overline{Pic^0 C}). $
\item[(0.9)]
$ I_j \in \overline{K} \ \forall j \geq 0, $
since
$ I_j = Hom_{\cO}(\cO^j,\cO) = (\cO^i)^v $.
In particular
$ \cM_Q \in \overline{Pic \ C}. \ $
\item[(0.10)]
$ \pi'_*(\overline{Pic^{-1} C'})\ \subseteq \ B(\overline{Pic^{0} C}) \ \ $
[J,(1.1)].
\item[(0.11)]
$ B(\overline{Pic^{0} C}) \ \subseteq \ \pi'_*(Pic^{-1} C)^{=} \ \ $ [J,(1.2)].
\item[(0.12)]
$ B(\overline{Pic^{0} C})
= \ \bigcup_i \pi'_*(\overline{Pic^{-i} C^i}) \
\Longleftrightarrow \ \cM = \cC \ \ $
[J,(1.3)].
\item[(0.13)] The following conditions are equivalent [J, $\S$2] : \\
(i) $ rk(\cM / \cM^2 + t \cC)  =  1 $ \\
(ii) $ B(\overline{Pic^{0} C}) = \pi'_{*} \overline{Pic^{-1} C'} $ \\
(iii) $ E(\cC, \delta) = \overline{K}  $ and
$ E(\cC', \delta -1) = \overline{K'}. \ $
\item[(0.14)] If
$ \ B(\overline{Pic^0 C}) = \pi'_*(Pic^{-1} C')^{=} \ $
then either $C$ is locally planar or
$ \ ord_t \cC = v_0 = 2 \delta -1 \ \ $ [J,\S 2]. $ \hfill{\Diamond} $
\end{enumerate}
\noindent{\underline{\bf  1: Decomposition.}}

In this section we give the basic decomposition of
$ \ (Pic^0 C)^{=} \ $
towards irreducible components. If $F$ is an $\cO$-module such that
$ \ \cO \ \subseteq \ F \ \subseteq \ \cO^{\sim}, \ $
and $ {\rm End}(F) $ does not contain $\cO'$, then
$ F' = uF, \ u \in \cO^{\sim{*}} $
for any
$ F' \approx F $
such that
$ \ \cC \ \subseteq \ F' \ \subseteq \ \cO^{\sim}. $
In particular
$ rk(F'/ \cC) = rk(F / \cC) $
is fixed. For $r$ such that
$ \ v_0 - \delta -1 \ \leq r \ \leq \delta -1, \ $
we define
\footnote{
$ M_{\tau} \ $
of [PS] is analogous to this.}
\smallskip\\
{\bf Definition (1.0):}
$ \ Pic^0_r C = \{ \cF \in (Pic^0 C)^{=} \mid rk(\cF / \cF \cC) = r+1 \}. \ $
Equivalently,
$ Pic^0_r C =
\{ \cF \in (Pic^0 C)^{=} \mid \ \exists F \approx \cF_Q \
{\rm such \ that} \
\cO \subseteq F \subseteq \cO^{\sim},
{\rm and} \
rk(F / \cC) = r+1 \}.$
\medskip\\
{\bf Lemma 1.1:} From definition of
$ Pic^0_r C, $ \\
(i) $ \ Pic^0_r C \ \bigcap \ Pic^0_{r'} C \ $
is empty for
$ r \neq r'. $ \\
(ii)
$ \ (Pic^0 C)^{=} \ =
\ \bigcup \ \{ Pic^0_r C \mid \ r =
\gamma-1, \cdots, v_0-2 \}. $
\footnote{[PS,Remark 3] is analogous to this.} \\
(iii)
$ \ \overline{Pic^0_r C} \ \subseteq \ Orb_{Pic^0 C} E(\cC, \delta-r) \
= \ \{ \cF \in (Pic^0 C)^{=} \mid \\ \exists \ F \ \approx \ \cF_Q \
{\rm \ such \ that} \ \cC \ \subseteq \ F \ \subseteq \ \cO^{\sim} \ {\rm and}
\ rk(F/\cC) = r+1 \}. \hfill{\parallel}$

We would like to know if and when
$ \ \overline{Pic^0_r C} \ $
is a component of
$ \ (Pic^0 C)^{=}, \ $
Generally
$ \overline{Pic^0_r C} $
is not irreducible.  However, since for all $r$ in this range
$ \ \exists \cF \in Pic^0_r C \ $
such that
$ End(\cF) $
does not contain
$ \cO_{C'}, $
\footnote{ (1.2), (1.3) (1.5) answer questions from Verdiere in private
talk (Pondicherry, Jan. '89) who pointed out that although these results
follow from the rest of this section, they deserve to be mentioned
independantly.}
\medskip\\
{\bf Lemma 1.2:}
$ \ \overline{Pic^0_r C} \ $
contains an irreducible component of
$ \ (Pic^0 C)^{=} \ $
for $r$ such that
$\gamma-2 < r < \delta. \hfill{\parallel}$
\medskip\\
Therefore
\medskip\\
{\bf Theorem 1.3:}  The minimum number of components of
$ \ Pic^0 C \ $
for a general unibranch singularity is
$ 2 \delta - v_0+1 $,
and these are the components containing an open subvariety
$ \ \{ \cF \mid End(\cF) \ {\rm does \ not \ contain} \ \cO_{C'} \}.
\hfill{\parallel}$

This bound is strict for some types of curves such as one with a non
Gorenstein singularity of embedding dimension greater than two, and
becomes an equality for locally planar or
$ \cM = \cC $
type
singularities.
\medskip\\
We generalize [R,1.2]: \ \ let
$ (-)^* $
denote
$ \ Hom(-, \cO_{C^{\sim}}) \ $
and let
$ K_r = ker(\pi^*_r) $
where
$ \ \pi^*_r : \ Pic^0_r \ \longrightarrow \ Pic^0 C^{\sim} \ $
is given by
$ \cF \ \mapsto \ \pi^* (\cF)^{**}. $
\medskip\\
{\bf Theorem 1.4:}
$ \ \cF \ \in \ \overline{Pic^0_r C} \
\Longleftrightarrow \ \exists \ \cL \ \in \ Pic^0 C \ $
such that
$ \ \cF \otimes \ \cL \ \in \ \overline{K}_r. \ $

\noindent{\bf Proof:}
Let $\cG$ be an
$ \ \cO_{C \chi Spec \ k[[t]]} $
-module such that
$ \ (\cG / t \cG) \ \approx \ \cF \ $
and $\cG$ restricted to the generic point of
$ \ Spec k[[t]] \ $
is generated by $r+1$ elements modulo
$ \cC $.
We have a diagram
$$
\begin{array}{lllll}
h    &   :  &   Spec \ k[[t]]  &           \ST           &  Pic^0 C,   \\[5mm]
h_0  &   :  &   Spec \ k(t)    &           \ST           & Pic^0_r C   \\
 {}  &  {}  &       {}         &           {}            & \DT \pi^*_r  \\
 {}  &  {}  &      Pic^0 C     &  \stackrel{\pi^*}{\ST}  & Pic^0 C^{\sim}
\end{array}
$$
with the generic point of
$ {\rm im}(h) $
given  by
$ h_0. $
Since
$ \ Pic^0 C^{\sim} \ $ is complete,
$ \pi^{*}_{r} o h_0 $
lifts to a map
$ \ p_0 : Spec \ k[[t]] \longrightarrow Pic^0 C^{\sim}, \ $
and since
$ \pi^* $
is smooth and surjective,
$ \ {\rm im}(p_0) \ $
lifts to a curve in
$ Pic^0 C $,
defining a line bundle
$ \cK^{-1} $
on
$ \ C \chi Spec \ k[[t]]. \ $
Let
$ \ \cP \ = \ \cK \otimes \cG. \ $
Then
$ \ (\cP / t \cP) \ \approx \ \cL \otimes \cF \ $
where
$ \ \cL \ = (\cK / t \cK) \ $
and $\cP$  restricted to the generic point of
$ \ Spec \ k[[t]] \ $
is in
$ K_r $.
This proves one direction, and since the
$ Pic^0 C $
-orbit of $K_r$
is contained in
$ \overline{Pic^0_r C} \ $
by definition, it proves the result.
$\hfill{\parallel}$
\medskip \\
Finally, since
$ \ \pi_* (Pic^{- \delta} C^{\sim}) \subseteq B(\overline{Pic^0_r C}) \
\ \forall r \ \in \ \{ 0, \cdots, v_0-2 \}, \ $
we have
\medskip\\
{\bf Lemma 1.5:}
$ \ (Pic^0 C)^{=} \ $
is connected.
\footnote{This is [PS, Th. 2].}
$ \hfill{\parallel}$
\newpage
\noindent{\underline{\bf 2: $ \cM = \cC:$}}

If
\footnote{The curve type
$ \cM = \cC $
of our notation is the type
$ \Gamma_{m,1,1} $
of [PS].}
$ \ \cM = \cC \ $
then
$ \ \overline{Pic^0_r C} \ $
is a component of
$ \ (Pic^0 C)^{=} \ $ for
$ r = 0, \cdots, v_0 -2, \ $
and we look at the intersections of components.
(From this section on readers are refered to [J], [PS] for proofs.)
\medskip\\
{\bf Theorem 2.1:}  If
$ \ \cM = \cC \ $
then
$ \ \overline{K}_r = E(\cC, \delta-r-1). \hfill{\parallel}$
\smallskip\\
Consequently if
$ \ \cM = \cC \ $
then
$ \ \overline{Pic^0_r C} = \ {\rm Orb}_{Pic^0 C} E(\cC, \delta -r-1).$
\medskip\\
{\bf Theorem 2.3:} If
$ \ \cM = \cC \ $
then
$ \ \overline{Pic^0_r C} \ $
is an irreducible component of
$ \ (Pic^0 C)^{=} \ $
for
$ \ r = 0, \cdots, v_0 - 2. \hfill{\parallel}$
\medskip\\
{\bf Lemma 2.4:} If
$ \ \cM = \cC \ $
then \\
(i)
$ \ B(\overline{Pic^0_r C}) = \pi' _{*} (\overline{Pic^{-1}_r C'}). $ \\
(ii)
$ \ Pic^0_r C \ \bigcap \ \overline{Pic^0_{r+1} C} \ = \ \phi. $ \\
(iii)
$ \ \overline{Pic^0_r C} \ \bigcap \ Pic^0_{r + j} C \
= \ \pi^{j}_{*} (Pic^{-j}_r C_j) \ $
for
$ r = 1, \cdots, v_0 -2 $
and
$ j = 1, \cdots, v_0 -2-r. \hfill{\parallel}$
\medskip\\
We  summarize the section in a couple of graphs.  Arrows indicate that the
object at source of arrow is contained in the object at the tip of the
arrow.

$$
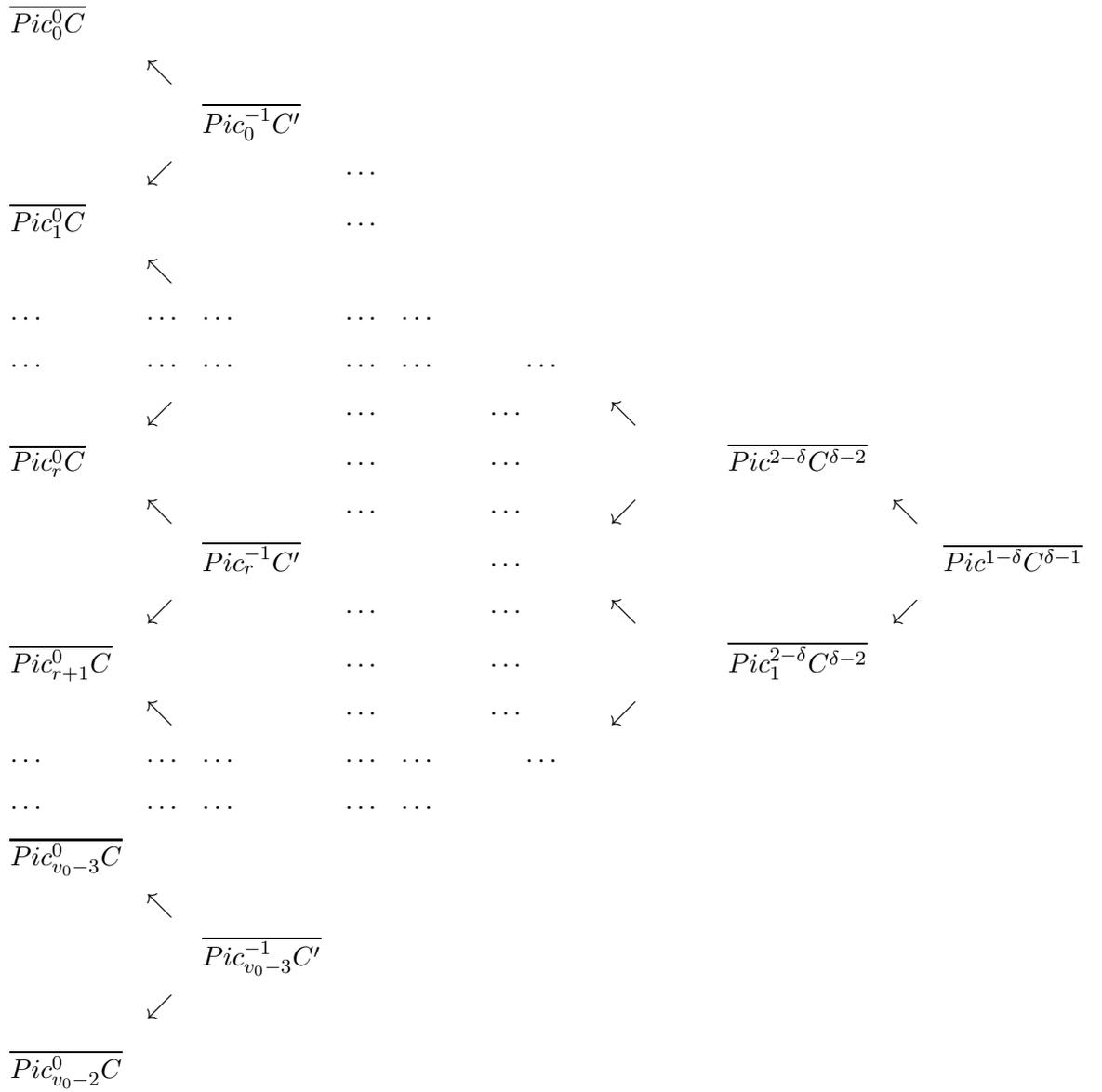
\begin{array}{llllllll}
Pic^0_0 C  &     {}     & {} & {} & \hspace{1in} {} & {} & {} & {} \\[2mm]
   {}      & $$\nwarrow$$ & {} & {} & \hspace{1in} {} & {} & {} & {} \\[2mm]
   {}      &     {}     & Pic^{-1}_0 C' &{}& \hspace{1in} {}&{}&{}&{} \\[2mm]
   {}      &     {}     & {} &{} & {} & {} & {} & {}\\[2mm]
Pic^0_1 C  &     {}     & {} & {} & \hspace{1in} {} & {} & {} & {} \\[2mm]
   {}      & $$\nwarrow$$ & {} & {} & \hspace{1in} {} & {} & {} & {} \\[2mm]
   {}      &     {}     & Pic^{-1}_1 C' &{}& \cdots
                                           \hspace{1.5in} {}&{}&{}&{} \\[2mm]
\cdots & \cdots & \cdots & \cdots & \cdots \hspace{1.5in} {}&{}&{}&{} \\[2mm]
\cdots & \cdots & \cdots & \cdots & \cdots \hspace{1.5in} {}&{}&{}&{} \\[2mm]
\cdots & \cdots & \cdots & \cdots & \cdots \hspace{.5in} \cdots \hspace{.5in}
                                               $$\nwarrow$$ {}&{}&{}&{}\\[2mm]
Pic^0_r C  &{}& \cdots & \cdots & \hspace{.5in} \cdots \hspace{.5in} {}&
                                     Pic^{2-\delta} C^{\delta-2} &{}&{}\\[2mm]
{}& $$\nwarrow$$ &{}& \cdots & \hspace{.5in} \cdots
                                 \hspace{.5in} {}&{}& $$\nwarrow$$ &{}\\[2mm]
{} & {} & Pic^{-1}_r C' &{}& \hspace{.5in} \cdots \hspace{.5in}
                     $$\nwarrow$$ {}&{}&{}& Pic^{1-\delta} C^{\delta-1}\\[2mm]
   {}      &     {}     & {} & {} & {} & {} & {} & {}\\[2mm]
Pic^0_{r+1} C &{}& \cdots & \cdots & \hspace{.5in} \cdots \hspace{.5in} {}&
                                  Pic^{2-\delta}_1 C^{\delta-2} &{}&{}\\ [2mm]
{} & $$\nwarrow$$ &{}&{}& \hspace{.5in} \cdots \hspace{.5in} {}&{}&{}&{}\\[2mm]
\cdots & \cdots & \cdots & \cdots & \cdots \hspace{.5in} \cdots
                                            \hspace{.5in} {}&{}&{}&{}\\[2mm]
\cdots & \cdots & \cdots & \cdots & \cdots \hspace{1.5in} {}&{}&{}&{} \\[2mm]
\cdots & \cdots & \cdots & \cdots & \cdots \hspace{1.5in} {}&{}&{}&{} \\[2mm]
Pic^0_{v_0 -3} C &{}&{}&{}& \hspace{1in} {}&{}&{}&{} \\[2mm]
{} & $$\nwarrow$$ & {} & {} & \hspace{1in} {}&{}&{}&{} \\[2mm]
{}&{}&  Pic^{-1}_{v_0-3} C &{}& \hspace{1in} {}&{}&{}&{}\\ [2mm]
   {}      &     {}     & {} & {} & {} & {} & {} & {}\\[2mm]
Pic^0_{v_0-2} C  &{}&{}&{}& \hspace{1in} {}&{}&{}&{}
\end{array}
$$

$$
\begin{array}{llllllll}
\overline{Pic^0_0 C} &{}&{}&{}& \hspace{1in} {}&{}&{}&{} \\ [2mm]
{}& $$\nwarrow$$ &{}&{}& \hspace{1in} {}&{}&{}&{} \\[2mm]
{}&{}& \overline{Pic^{-1}_0 C'} &{}& \hspace{1in} {}&{}&{}&{} \\[2mm]
{}& $$\swarrow$$ &{}& \cdots & \hspace{1in} {}&{}&{}&{} \\[2mm]
\overline{Pic^0_1 C} &{}&{}& \cdots & \hspace{1in} {}&{}&{}&{} \\[2mm]
{}& $$\nwarrow$$ &{}&{}& \hspace{1in} {}&{}&{}&{} \\[2mm]

\cdots & \cdots & \cdots & \cdots & \cdots \hspace{1.5in} {}&{}&{}&{} \\[2mm]
\cdots & \cdots & \cdots & \cdots & \cdots \hspace{.5in} \cdots
                                           \hspace{.5in} {}&{}&{}&{} \\[2mm]
{} & $$\swarrow$$ & {} &  \cdots  & \hspace{.5in} \cdots \hspace{.5in}
                                             $$\nwarrow$$ {}&{}&{}&{} \\[2mm]

\overline{Pic^0_r C} &{}&{}& \cdots & \hspace{.5in} \cdots
     \hspace{1in} {} & \overline{Pic^{2-\delta} C^{\delta-2}} &{}&{}\\[2mm]
{}& $$\nwarrow$$ &{}& \cdots & \hspace{.5in} \cdots \hspace{.5in}
                              $$\swarrow$$ {}&{}& $$\nwarrow$$ &{}\\[2mm]
{}&{}& \overline{Pic^{-1}_r C'} &{}& \hspace{.5in} \cdots \hspace{.5in}
                 {}&{}&{}& \overline{Pic^{1-\delta} C^{\delta-1}} \\[2mm]
{}& $$\swarrow$$ &{}& \cdots & \hspace{.5in} \cdots \hspace{.5in}
                                   $$\nwarrow$$ {}&{}& $$\swarrow$$ &{}\\[2mm]
\overline{Pic^0_{r+1} C} &{}&{}& \cdots & \hspace{.5in} \cdots
     \hspace{.5in} {} & \overline{Pic^{2-\delta}_1 C^{\delta-2}} &{}&{}\\[2mm]
{}& $$\nwarrow$$ &{}& \cdots & \hspace{.5in} \cdots \hspace{.5in}
                                     $$\swarrow$$ {}&{}&{}&{} \\[2mm]
\cdots & \cdots & \cdots & \cdots & \cdots \hspace{.5in} \cdots
                                             \hspace{.5in} {}&{}&{}&{}\\[2mm]
\cdots & \cdots & \cdots & \cdots & \cdots \hspace{1.5in} {}&{}&{}&{} \\[2mm]

\overline{Pic^0_{v_0 -3} C} &{}&{}&{}& \hspace{1in} {}&{}&{}&{} \\[2mm]
{}& $$\nwarrow$$ &{}&{}& \hspace{1in} {}&{}&{}&{} \\[2mm]
{}&{}& \overline{Pic^{-1}_{v_0-3} C'} &{}& \hspace{1in} {}&{}&{}&{} \\ [2mm]
{}& $$\swarrow$$ &{}&{}& \hspace{1in} {}&{}&{}&{} \\[4mm]
\overline{Pic^0_{v_0-2} C} &{}&{}&{}& \hspace{1in}{}&{}&{}&{}
\end{array}
$$
\newpage
\underline{\bf 3: $ \cM' = \cC', \ C $ Gorenstein:}

$\cO$  is Gorenstein if and only if
every torsion free $\cO$-module is reflexive, i.e, if and only if
$ \ \Gamma = \{ j \mid \ v_0-1-j \not \in \Gamma \}. $
So for all
$ \ \cF \in Pic^0 C, $
we have
$ \ End(\cF) = \cO_C \Longleftrightarrow
\cF_Q \approx \cO, \Longleftrightarrow \cF \in Pic^0 C. \ $
In particular therefore the components of
$ \ (Pic^0 C)^{=} \ $
other than
$ \ \overline{Pic^0 C}, \ $
if any, are components of
$ \ \pi'_* (Pic^{-1} C')^{=}. \ $
If $C$ is locally planar, we know that
$\ \overline{Pic^0 C} \ = (Pic^0 C)^{=} \ $
and
$ \ \pi'_* (Pic^{-1} C')^{=} \ $
is all of
$ \ B(\overline{Pic^0 C}) $,
contributing no components to
$ \ (Pic^0 C)^{=}. \ $
On the other hand, if
$ \ \cM' = \cC' \ $
and if $C$ is not locally planar
\footnote{ This is
$ \Gamma_{m,m-1} $
of [PS].}
then
$ \Gamma = (k_1, \cdots, 2k_1-2){\bf N}, \ k_1 \geq 3. $
\medskip\\
{\bf Theorem 3.1:}  Let  $C$  be Gorenstein, and not locally planar.  If
$ \ \cM' = \cC', \ $
then \\
(i) $ \ (Pic^0 C)^{=} \ = Pic^0 C \ \bigcup \ \pi'_* (Pic^{-1} C')^{=} \ , \ $
and components of
$ \ (Pic^0 C)^{=}, \ $
other than
$ \ \overline{Pic^0 C}, \ $
are precisely
$ \{ \pi'_* (\overline{Pic^{-1}_r C'}) \mid r = 1, \cdots, \delta -2 \}. $ \\
(ii) If
$ \ \cF \in \overline{Pic^0 C}, \ $
and if
$ \ F \approx \cF_Q \ $
is such that
$ \ F \in E(\cC, \delta), \ $
then
$ \ \forall \ r \geq 1, \ $
we have
$ ord_t F = r \
\Longrightarrow \cF \in \pi'_{*} (\overline{Pic^{-1}_{r-1} C'}). \ $
Conversely, for all
$ \ \cF \in \pi'_{*} (\overline{Pic^{-1}_{r-1} C'}), \ \
\cF \in \overline{Pic^0 C} \
\Longrightarrow \exists F \approx \cO_Q \ $
such that
$ \ \cC \subseteq F \subseteq \cO^{\sim}, \ $
with
$ ord_t F = r, \ F_1 \subseteq t^{k_1} \cO^{\sim} \ $
and
$ \ rk(F_1 / t^{k_1+r}\cO^{\sim}) = r-1. \hfill{\parallel}$
\medskip\\
\underline{\bf Example 1:} Let $\cO$  be generated by monomials
in some local parameter $t$, and let
$ F = (t^2, t^{k_1}) \cO + \cC $.
The only possible deformation of $F$ is given by either
$ (t^2 + \beta) \cO $
when $k_1$ is even or
$ (t^2 + \beta t + \beta^2) \cO $
when $k_1$ is divisible by 3, and there is no deformation if
$ k_1 \equiv \pm 1$ modulo 6.
$\hfill{\diamond}$
\medskip\\
So the reverse implications of (ii) are not true. On the other hand
\newpage
\underline{\bf Example 2:}  For
$ \cO = k[[t^4, t^5, t^6]] $
we have
$ Filt(\cC, \delta) \ = \ \overline{K} $:
if
$ ord_t F = 1 $
then
$ F = t\cO'; $
if
$ ord_t F = 2 $
and
$ F \neq t^2 \cO^{''} $
then
$F = (t^2 , ut^4) + \cC$
for
$u \in \cO^{\sim{*}}$,
and
$\partial_{\beta} = t^2 + \beta u$
will do; if
$ord_t F = 3$
then $F$ is in the closure of the Picard orbit of
$(t^3, t^4, t^5)+\cC = t^{-1} \cM$, and
$ ord_t F = 4 $
if and only if
$ F = \cO^{\sim} $.
In any case $F$  has the desired deformation.
$\hfill{\diamond}$
\medskip\\
\underline{\bf 4: $ Rk(\cM / \cC) = 1, \ \cM' = \cC' $:}

 From (0.13), if
$ rk(\cM / \cC) = 1 $,
then
$ Filt(\cC, \delta) = \overline{K}. $
By definition each
$ \overline{Pic^0_r C} $
contains at least one component of
$ (Pic^0 C)^{=} \ $
for
$ j \in \{ 1 , \cdots , k_1 -1 \}. $
Since
$ (Pic^0 C)^{=} \ = \ \bigcup \ \{ \ \overline{Pic^0_r C} \
\mid r = 1, \cdots, k_1 -1 \}, \ $
we need to find components of
$ \overline{Pic^0_r C} \ $
and to see when these are actually components of
$ (Pic^0 C)^{=}. \ $
If
$ Rk(\cM / \cC) = 1 $
and
$ \cM' = \cC' $
\footnote{ If $C$ is a monomial singularity curve then this is
$ \Gamma_{m,1,2} $
of [PS].} \
then we have
$ \ \overline{Pic^0_0 C}  = \phi, \ Pic^0_1 C  = Pic^0 C,
\ Pic^0_2 C  = \{ \cF \mid \cF_Q = \cO + f_1 \cO, \
ord_t f_1 \geq 2 \}, $
and
$ Pic^0_{k_1} C = {\rm orb}_{Pic^0 C} \cW^0_C. \ $
Let
$ \cF \in Pic^0_r C \ $
and
$ F \approx \cF_Q \ $
such that
$ \ u \cO \subseteq F \subseteq \cO^{\sim} \ \
{\rm and} \ u \in \cO^{\sim{*}}, \ $
so that
$ rk(F / \cC) = r+1. \ $
Let
$ \{ f_0, f_1, \cdots, f_r \} $
be a basis of
$ \ F, \ f_i = \sum \{a_{ij} t^j \mid j = 0,
\cdots, k_1 +1, a_{ij} \in k \}, \ $
and, say,
$a_{00} \neq 0$.  Then
we must have
$ \ rk[((f_0, \cdots f_r) \cO + \cC) / \cC] \ = r+1, $
$$
\noindent{\begin{tabular}{lll}
$\Longrightarrow$
& rank of $M_F$, the coefficient matrix of $(F/\cC)$, is $r + 1$.\\
$\Longrightarrow$
& every $ (r+2) \times (r+2) $ minor of $M_F$  must vanish.  \\
\end{tabular}}
$$
Therefore
$ \ \overline{Pic^0_r C}  =
\overline{P_{r,1}} \ \bigcup \ \overline{P_{r,2}}, \ $
where
$$
\begin{tabular}{lll}
$ P_{r,1} $ &  = &
$ \{ \cF \in Pic^0_r C \mid \exists F \approx \cF_Q \ {\rm such \ that} $ \\
      {}    & {} &
$ F \in E(\cC, \delta - r-1), \ \cO \subseteq F, \ {\rm and} \ ord_tF_1 = 1 \}$
\end{tabular}
$$
and
$$
\begin{tabular}{lll}
$ P_{r,2} $ &  = &
$ \{ \cF \in Pic^0_r C \mid \exists F \approx \cF_Q \ {\rm such \ that} $ \\
    {}      & {} &
$ F \in E(\cC, \delta-r-1), \
\cO \subseteq F, \ {\rm and} \ ord_t F_1 \geq 2 \}.$
\end{tabular}
$$
It remains to see whether or not
$ \overline{P_{r,i}} \ $
are components of
$ Pic^0 C \ $
for
$ r = 1, \cdots, k_1 -1 $.
Since
$ P_{r,1} = \pi'_{*} Pic^{-1}_{r-2} C' \ $
for
$ r = 2, \cdots, k_1, \ \ \overline{P_{r,1}} \ $
is not contained in
$ P_{r',1} $
for
$ r \neq r' $.
Hence proving that
$ P_{r,i} $
are, indeed, components of
$ \ (Pic^0 C)^{=}, \ $
is now reduced to the question of looking at
$ B(\overline{P_{r,2}}) $.
\smallskip\\
{\bf Theorem 4.1:}
$ \ \cF \in \overline{P_{r,2}} $
if and only if the following conditions hold:\\
(i) $ \ \exists F \ $
such that
$ \ \cF_Q \approx F, \ \cC \subseteq F
\subseteq \cO^{\sim}, \ rk(F / \cC) = r+1, \ $
and \\
(ii)  either
$ ord_t F = 1 \ {\rm and} \ F \supseteq t^{k_1+1}
\cO^{\sim}, \ $
or
$ \ ord_t F \geq 2 $,
where \\
$ \ F \supseteq F_1 \supseteq \cdots \supseteq F_r \supseteq \cC \ $
is a filtration with the following properties:
$ \ rk(F_i / F_{i+1}) = 1, \ ord_tF_{i+1} \ > \ ord_t F_i, $
and
$ \ F_r \subseteq t^{k_1} \cO^{\sim}. $
$\hfill{\parallel} $
\smallskip \\
{\bf Corollary 4.1.1:}
$ \ P_{i,1} \subseteq \ \overline{P_{i-1, 2}} \ $
for all
$ i < k_1. \hfill{\parallel}$
\smallskip\\
{\bf Corollary 4.1.2:}
$ \ \cF \in \overline{P_{r,2}} \ \bigcap \ \overline{P_{r+1,2}} \
\Longleftrightarrow \ \exists F \approx \cF_Q \ {\rm such \ that} \
\cC \subseteq F \subseteq t\cO^{\sim}, \ rk(F/\cC) = r+1 \
{\rm and} \ F_1 \subseteq t^3 \cO^{\sim}. \hfill{\parallel}$
\smallskip\\
For
$ i = 1, \cdots, k_1 -1, \ \exists \cF \in P_{i,2} \
{\rm such \ that} \
\cF_Q \approx \sum \{ t^j \cF \mid j = 0,2, \cdots, i \} $.
Since
$ \cF  \not \in P_{i',2} \ \forall \ i \neq i' $,
this proves
\smallskip\\
{\bf Theorem 4.2:}
$ (Pic^0 C)^{=} \ = \
\bigcup \ \{ P_{i,2} \mid i = i, \cdots, k_1 -1 \} $
and
$ \ \overline{P_{i,2}} \ $
is a component of
$ (Pic^0 C)^{=} \ $
for
$ i = 1, \cdots, k_1-1.  \hfill{\parallel} $
\medskip\\
\begin{center}
\underline{\bf References}
\end{center}
\begin{enumerate}
\item[{[AIK]}] A. Altman, A. Iarrobino and S. Kleiman,
Irreducibility of the compactified Jacobian, in ``Proceedings,
Nordic Summer School NAVF, Oslo, Aug. 5-25, 1976,'' Noordhoff,
Groningen, 1977.
\item[{[AK1]}] A. Altman and S. Kleiman, Compactifying the
Picard Scheme, Advances In Mathematics 35, 50-112 (1980).
\item[{[AK]}] A. Altman and S. Kleiman, Compactifying the Picard
Scheme, II, Amer. J. Math. 101 (1979) (``Zariski volume''),
10-41.
\item[{[D]}] C. D'Souza, ``Compactification of Generalized
Jacobians,'' Thesis, Tata Institute, Bombay, 1973.
\item[{[HK]}] Herzog and Kunz, Die Wertehalbgruppe eines Lokalen
Ringes der Dimension I, Springer Verlag, 1971.
\item[{[J]}] $\underline{~~~~~~~}$ , Thesis, Brandeis University, 1986.
\item[{[J']}] $\underline{~~~~~~~}$, Boundary of Picard Group of a
Singular Curve, preprint.
\item[{[K]}] H. Kleppe, Picard scheme of a curve and its
compactification.  Thesis, M.I.T., 1981.
\item[{[MM]}] A. Mayor and D. Mumford, Further Comments on
Boundary Points, Amer. Math. Soc. Summer Institute, Woods Hole,
Mass., 1964.
\item[{[PS]}] G.Pfister and J.H.M. Steenbrink, Reduced Hilbert Scheme for
Irreducible Curve Singularities, to appear, Journal of Pure and Applied
Algebra.
\item[{[R]}] C. Rego, The Compactified Jacobian, Ann. Scient.
Ec. Norm. Sup., 4$^4$ serie, t. 13, 1980, p. 211 a 224.
\item[{[SK]}] S. Kleiman: The structure of the compactified
Jacobian : a review and an announcement; seminari di geometria,
1982-83; Universita delgi Studi di Bologna, Dipartmento
Mathematica, 1984.
\end{enumerate}
\end{document}